\documentclass[12pt,oneside,a4paper]{article}

\usepackage{subfigure}
\usepackage{latexsym}
\usepackage{graphicx}

\title{State space description of national economies:\\
       the V4 countries}

\author{Ivo Petras, Igor Podlubny \\[1ex]
        Department of Applied Informatics and Process Control, \\
        BERG Faculty, Technical University of Kosice, \\
        B. Nemcovej 3, 04200 Kosice, Slovak Republic\\[1ex]
        e-mail: ivo.petras@tuke.sk, igor.podlubny@tuke.sk}

\date{\today}

\begin{document}
\maketitle

\begin{abstract}

We present a new approach to description of
national economies. For this we use the state space viewpoint,
which is used mostly in the theory of dynamical systems
and in the control theory. Gross domestic product,
inflation, and unemployment rates are taken as state
variables. We demonstrate that for the considered
period of time the phase trajectory of each of the V4 countries
(Slovak Republic, Czech Republic, Hungary, and Poland)
lies approximately in one plane, so that the economic
development of each country can be assocated with
a corresponding plane in the state space.
The suggested approach opens a way to a new set
of economic indicators (for example, normal vectors
of national economies, various plane slopes,
2D angles between the planes corresponding to
different economies, etc.).

The tool used for computations is orthogonal regression
(alias orthogonal distance regression,
alias total least squares method), and we also
give general arguments for using orthogonal regression
instead of the classical regression
based on the least squares method.

A MATLAB routine for fitting 3D data to lines
and planes in 3D is provided.

\vspace*{2ex}


\end{abstract}

\section{Introduction}

In this article we consider two different versions
of the least squares method -- the classical one
and the so called orthogonal regression.
We incline to using the orthogonal regression,
and give several arguments in favour of this.
The main reason is the fact that orthogonal regression
is a suitable tool for fitting lines and surfaces
in multidimentional space, while even the use
of classical regression in 2D is not so natural.

The article is organized as follows.

First, we briefly recall the history of the classical
least squares method to remind that the justification
of its applications is unclear from its very beginning.

Second, we take a critical look at the traditional interpretation
of the classical least squares method and conclude
that it is absolutely artificial and meaningless.
Changing a viewpoint from ``squares'' to ``circles'',
we find a natural geometric interpretation for
quantities to be minimized, and this change
leads to the method of orthogonal regression.

After that, we provide several general arguments
in favour of using orthogonal regression.

We provide an easy example -- ``The Flight of a Bumblebee'' --
in order to illustrate how the orthogonal regression method
can help in describing dependences in multidimentional
space.

Then we apply the orthogonal regression method to describing
the national ecomonies of the countries of the V4 group
in the state space of three variables -- gross domestic product,
inflation, and unemployment. The development of the national
economies under study is described by their trajectories in
the chosen state space. We discovered that for each particular
country the trajectory of its national economy lies
approximately in one plane, so the normal vector of this
plane and one of its points (for example, the one coinciding
with the data centroid) can be associated with a particular
economy.

In the conclusion we discuss some new possibilities that
our approach opens in the field of description of
complex dynamic economic systems such as national economies.

\section{The history of the least squares method}

One can hardly find another topic, which is so famous
and which has such misterious and unclear history
as the least squares method.

The controversal circumstances around the
the appearance of the method of least squares were
investigated by many authors.
According to Celmins \cite{Celmins-1998},
it might be summarized as follows.

In 1805 Legendre published ``Nouvelles methodes pour la determination des
cometes'', in which he introduced the method of least squares
and gave it this famous name.

In 1809, Gauss published the book ``Theoria motus corporum coelestium in
sectionibus conicis solem ambientium'' \cite{Gauss-1809,Gauss-1963}, where he discussed
the method of least squares and, mentioning Legendre's work, stated that he
himself had used the method since 1795.

Legendre felt offended by Gauss's statement.
In a letter to Gauss  about his new book Legendre
wrote that claims of priority should not be made without proof by previous
publications. Gauss did not have such a publication,
but he stated that he was convinced that the
idea of least squares method is so simple that many people must have used
the method even before him.

Later Gauss tried to prove his claim but had only little success.
His own computational notes were -- as he said -- lost.
His colleagues apparently did not remember discussions with him
or did not want to be involved in the controversy.
As Celmins mentions \cite{Celmins-1998},
only the astronomer Olbers included in a paper
in 1816 a footnote asserting that Gauss had shown him the method of least
squares in 1802, and Bessel published a similar note in a report in 1832.

Later Gauss gave up the search but did not retract his claim. In 1820 Legendre
published a supplement to his 1805 memoir with an appendix where he publicly
attacked Gauss's claims of priority. The controversy continued, and in 1831
H.~C.~Schumacher (cf. Celmins \cite{Celmins-1998})
wrote to Gauss about a publication published in 1799 that contained data and
adjustment results by Gauss. Schumacher suggested repeating the calculations and
thereby demonstrating that the method of least squares was indeed used by Gauss
in 1799. Gauss's answer was that he was well aware of the data but would not
permit a recalculation, since such attempts would only
suggest that he could not be trusted.

After describing the above circumstances, Celmins tried to repeat
Gauss's calculations.
Celmins's main conclusion was that Gauss's results are not
consistent with any obvious and reasonable adjustment of
observational errors nor with a least-square adjustment
in the parameter space.

\section[Classical LSM]{Classical LSM}

The well known least squares method is
a mathematical procedure for finding the best-fitting curve to a given set of
points by minimizing the sum of the squares of the offsets (``the residuals'') of
the points from the curve.

The classical least squares fitting consists in minimizing
the sum of the squares of the vertical deviations  of a set
of  data points
\begin{equation}\label{eq:LSM-criterion}
E = \sum_{i} [y_i - f(x_i, \alpha_1, \alpha_2, \ldots, \alpha_n)]^2
\end{equation}
\noindent
from a chosen function $f$.

%
%
%

\begin{figure}
\centering
\subfigure[Fitting data to a straight line]{\label{fig:fitting-data-to-line}\includegraphics[width=0.45\textwidth,angle=0]{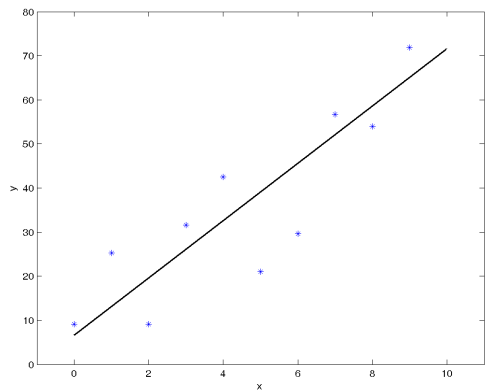}}
\hfill
\subfigure[Geometric interpretation(?) of the classical LSM]{\label{fig:fitting-by-LSM}
                   \includegraphics[width=0.45\textwidth,angle=0]{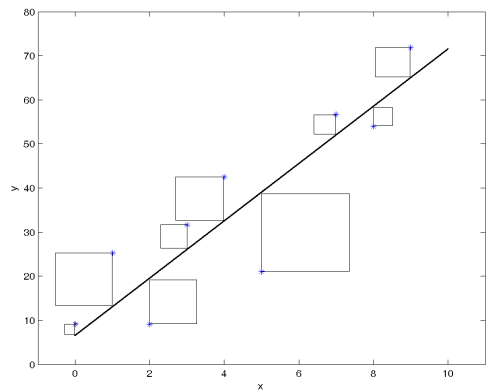}}
  \caption{Fitting by the classical least squares method (LSM)}
\end{figure}

For a simple illustration, let us recall the classical
linear regression problem, in which we have to fit
the set of data by a straight line in 2D plane.
This situation is shown in Fig.~\ref{fig:fitting-data-to-line}.

One may use the edge of a transparent ruler or a tightly stretched
black string to get a line which seems to fit the data points.
Mathematically, what we have to do
is to determine the parameters $k$ and $b$ of the
equation of a straight line:
\begin{equation}\label{eq:lin-regr-example}
y = k x + b
\end{equation}
where
$x$ -- the variable considered as independent;
$y$ -- the variable considered as dependent on $x$;
$k$, $b$ -- constants, often called parameters,
                to be determined so that the line fits the
		data optimally (in some sense...).

To find a way to calculate the parameters, let us return to the simplest case,
equation (\ref{eq:lin-regr-example}) and
Fig.~\ref{fig:fitting-data-to-line},
and let us agree to choose $k$ and $b$ so as to minimize the sum of the squares of the
errors. Fig. \ref{fig:fitting-by-LSM} shows these squares graphically.
We twist and push the line through these
points until the sum of the areas of the squares is the smallest.
The natural question is: ``Why the squares?
Why not just get the smallest sum of distances of the data points from the line?''
The only real explanation is that
it is easy to compute $k$ and $b$
to minimize the sum of squares off offsets of $y$ (vertical offsets),
but it is quite difficult to minimize
(using analytic derivations) the sum of distances of the data points from the line,
and it is really this ease that is responsible
for the generally accepted preference for the squares of offsets of $y$.

In fact, we must always keep in mind that the least squares
approach is basically a ``last resource'' tool that is used for
obtaining at least some mathematical model for a process under study,
when obtaining a better model by analytical derivations is
impossible or extremely time and/or effort consuming.
As such, it is based on several presumptions, which
are based only on some intuition or prejustice;
among them we must mention postulating the role
of variables (``independent/dependent'') and the type
of dependance between them (linear, polynomial, exponential,
logarithmic, etc.)

\section{The method of least circles?!} \label{sec:LCM}

Looking at the geometric ``interpretation'' of the least squares method
shown in Fig.~\ref{fig:fitting-by-LSM}, which so often appears
in numerous textbooks and lectures, we can conclude
that it is absolutely artificial and does not contain
any sign of mathematical beauty. The picture shown
in Fig.~\ref{fig:fitting-by-LSM} again provokes the
question: ``Why the squares?''

To change a viewpoint, let us note that the criterion
(\ref{eq:LSM-criterion}) can be painlessly replaced
with
\begin{equation}\label{eq:LSM-criterion-beta}
E = \beta \sum_i [y_i - f(x_i, \alpha_1, \alpha_2, \ldots, \alpha_n)]^2
\end{equation}

Indeed, multiplication by a non-zero number $\beta$ does not affect
the point of minimum. Only the minimum value of the criterion
function ($E$) will be multiplied by $\beta$ -- but this
value itself is not the subject of interest, since we look
for the values of $\alpha_1$, $\alpha_2$, ... $\alpha_n$.

Taking $\beta=\pi$, we obtain
\begin{equation}\label{eq:LSM-criterion-pi}
E = \sum_i \pi [y_i - f(x_i, \alpha_1, \alpha_2, \ldots, \alpha_n)]^2
\end{equation}

\begin{figure}
\centering
\noindent
\subfigure[The case of classical least squares fitting]{\label{fig:LCMVert}\includegraphics[width=0.45\textwidth,angle=0]{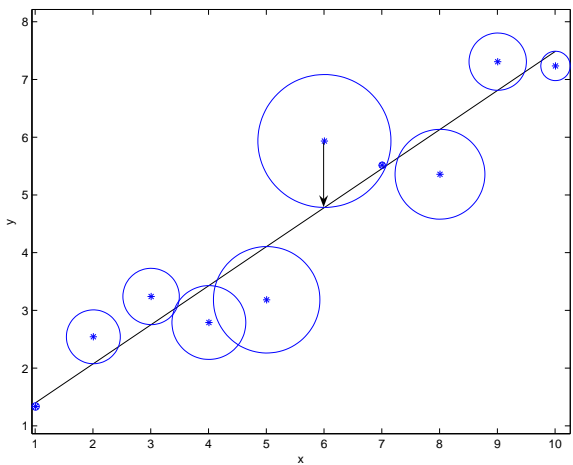}}
  \hfill
  \subfigure[The case of orthogonal distance fitting]{\label{fig:LCMOrth}
                   \includegraphics[width=0.45\textwidth,angle=0]{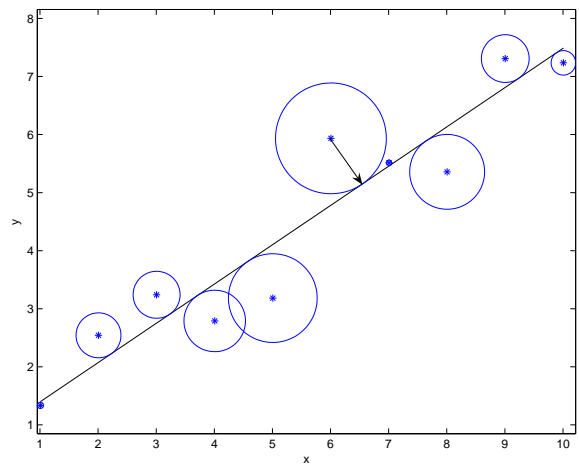}}
  \caption{``Least circles'' viewpoint}
\end{figure}

%

Geometrically, the formula (\ref{eq:LSM-criterion-pi}) means
the sum of areas of the circles shown in Fig.~\ref{fig:LCMVert}.
The radii of the circles in Fig.~\ref{fig:LCMVert} are
the vertical offsets of $y_i$ from the fitting line.
Each of those circles has two points of intersection
with the line. It is clear, that one cannot consider
this picture as elegant. Changing the radii slightly,
one can preserve $n$ pairs of intersection of the circles
and the line, so one can draw an infinite number of
pictures lloking similar to Fig.~\ref{fig:LCMVert}.
But Fig.~\ref{fig:LCMVert} is just a reformulation
of the standard geometric ``illustration'' of the least
squares method (recall Fig.~\ref{fig:fitting-by-LSM}).
Instead of the ``least squares method'' we now deal
with the ``least \emph{circles} method''. But the circles
shown in Fig.~\ref{fig:LCMVert} are clearly not the best.

However, the circles shown in Fig.~\ref{fig:LCMOrth}
are really optimal: the fitting line is a tangent line
to all circles.
The radii of the circles in
Fig.~\ref{fig:LCMOrth} are equal to minimal distances
between the points $(x_i, y_i)$ and the fitting line,
and this guarantees the unique picture.

The criterion to minimize in this case is
\begin{equation}\label{eq:LCM-criterion-pi}
E = \sum_i \pi \, \Bigl[ d\Bigl((x_i,y_i), f(x, \alpha_1, \alpha_2, \ldots, \alpha_n) \Bigr) \Bigr]^2,
\end{equation}
which is up to a constant multiplier $\pi$ the
formula known under the name of orthogonal regression
or total least squares \cite{de-Groen-intro-1996,Golub-van-Loan-1980,Huffel-Vandrwalle-1991-book,Nievergelt-review-1994}.
Here $d\Bigl((x_i,y_i), f\Bigr)$ denotes the distance between
the point $(x_i,y_i)$ and the fitting line $f$.

\section{Arguments for orthogonal regression approach}

There are numerous works and significant actvities
devoted to the theory of the total
least squares method \cite{Book-of-abstracts-TLSEIVM-2006}.
To support its wider applications,
we would like to list the following arguments
in favour of the orthogonal distance
(or total least squares) fitting.

\begin{enumerate}

\item
\emph{The shortest (orthogonal) distance is the most natural viewpoint on any fitting.}

The section~\ref{sec:LCM} and the
comparison of Fig.~\ref{fig:LCMVert} and Fig.~\ref{fig:LCMOrth}
give a clear visual evidence for this.

\item
\emph{The sum of orthogonal distances is invariant with respect to the choice of the system of coordinates.}

This is obvious, since the distance is independent of the choice of the coordinate system.
In addition, let us imagine that we replace rectangular coordinates
with non-rectangular. For example, in Fig.~\ref{fig:affine-coordinates}
an affine system with the angle of 60 degrees between the axes is shown.
In this picture the squares of $\Delta y_i$ have even less meaning than
in Fig.~\ref{fig:fitting-by-LSM}. However, the orthogonal distances,
which are also shown in Fig.~\ref{fig:affine-coordinates}
still have solid interpretation and can be used for fitting.

\begin{figure*}
\centering
\includegraphics[width=0.6\textwidth,angle=0]{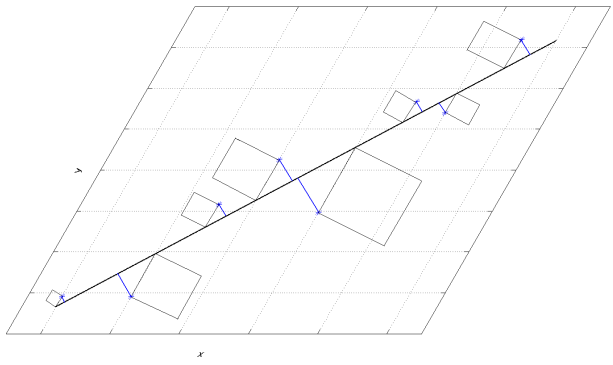}
\caption{The case of non-carthesian coordinates.}\label{fig:affine-coordinates}
\end{figure*}

\item
\emph{There are no conjugate regression lines, which appear after swapping $x$ and $y$, because in the case of orthogonal regression the fitting $y=f(x)$ gives exactly the same line as the fitting $x = f^{-1}(y)$.}

Suppose one wants to find the dependance between the height ($x$) and the weight ($y$) of people.
The dependance is presumed to be linear (a straight line) described by equation~(\ref{eq:lin-regr-example}).
After determining $k$ and $b$, this relaitonship can be used for estimating
the weight of a person of a given height.

However, the vewpoint can be inverted: for a given weight,
estimate the height of a person. If one already has the equation
(\ref{eq:lin-regr-example}), then the solution should normally be
\begin{equation}
x = \frac{1}{k} y - \frac{b}{k}.
\end{equation}

However, the classical regression approach leads to a different result
expressed by a conjugate regression line.

In Fig.~\ref{fig:LS-vs-TLS} a solution to the Nievergelt's
easily understandable example \cite{Nievergelt-review-1994} is shown.
To make the situation more obvious, we added one to the ordinates
in the Nievergelt's example,
so we used the data shown in Table~\ref{tab:Nievergelts-example-data}.
The classical least squares regression $y$ versus $x$ gives
the regression line $y =0.45 x + 3.2$. After swapping $x$ and $y$,
classical least squares regression gives
the conjugate regression line $x = 0.45 y + 1.75$.
However, in the case of the orthogonal distance regression we obtain
the same line $y = x+1$ ($x = y-1$) independently on the order of $x$ and $y$.

It is worth noting that all three lines run through the centroid,
and that the orthogonal distance regression line is located between
the ``scissors'' formed by the conjugate regression lines
obtained by the classical least squares regression.

\begin{table}
\caption{Data for the Nievergelt's example.}\label{tab:Nievergelts-example-data}
\begin{center}
\begin{tabular}{c|ccccc}
\hline
x	& 1	& 3	& 4	& 5	& 7 \\
y	& 4	& 2	& 6	& 8	& 5 \\
\hline
\end{tabular}
\end{center}
\end{table}

\begin{figure}
\centering
\includegraphics[width=0.6\textwidth,angle=0]{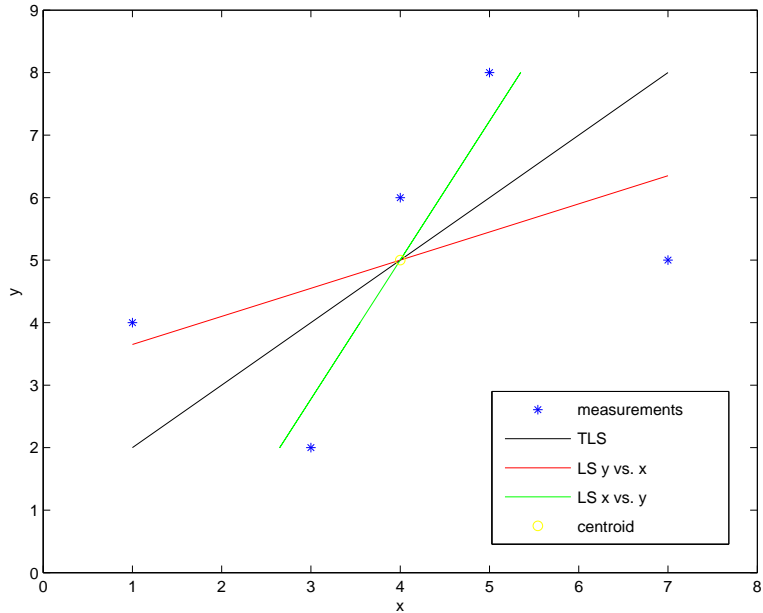}
\caption{Nievergelt's example: classical regression versus orthogonal regression.}\label{fig:LS-vs-TLS}
\end{figure}

\item
\emph{There are no problems with causality (normally, determination
of what is an independent variable and what is a dependent variable
is simply unclear or even impossible; this is always postulated).}

In all textbooks on statistics discussing the classical regression analysis
it is always underlined that that the choice of what is the ``independent'' variable
and what if the ``dependent'' variable is extremely important.
However, in many cases it is not so easy to make a decision what is what,
and justification of such a decision is based only on
some subjective judgement or on a prejustice.
We already mentioned the relationship
between the weight and the height of people.
It is obvious that any of these characteristics can be
taken as an independent one.

\item
\emph{Implementation of the orthogonal fitting
does not depend on the number of dimensions.}

To realize this, just recall the square of the distance between the two points
$P(x_1, x_2, \ldots, x_n)$ and $Q(y_1, y_2, \ldots, y_n)$
in  $n$-dimensional orthogonal coordinate system is
\begin{equation}
[d(P,Q)]^2 = \sum_{i=1}^{n} (x_i -y_i)^2
\end{equation}

\end{enumerate}

\section{State space description of national economies}

In everyday professional and non-professional communication
one can frequently hear the word ``state of economy''.
This expression can be given an exact meaning
by adopting the tools that are available in the
theory of dynamical systems and in automatic control.
Namely, we will use the technique called state space
description and phase trajectories.

Suppose we have a dynamical process (that is, the process
deveoping in time), which is characterized by three
quantities $x(t)$, $y(t)$, $z(t)$ (of course,
number or dimensions can be different, not necessarily 3).
Chosing $x$, $y$, and $z$ as the three coordinates,
we can assign a point $(x(t), y(t), z(t))$ to each
value of $t$ -- that is, we assign a point in a 3D space
to a state of the considered process at time $t$.
The variables $x$, $y$, $z$ are called the state variables.
The line formed by the points $(x(t), y(t), z(t))$
when $t$ takes on the values from a given interval
(usually $[0, T]$ for some finite $T$, or $[0, \infty)$)
is called the phase trajectory of the process.

As state variables for this study, we selected those
which are standard: gross domestic product (GDP),
inflation, and unemployment. In the subsequent sections,
we demonstrate the advantages of the state space
description of national economies and possible ways
for further developments and scientific investigations.

\section{Example: ``The Flight of the Bumblebee'' -- fitting 3D data by a straight line in 3D}

Before continuing with real data for the countries of the V4 group,
let us consider a simple general example.
Suppose a bumblebee is flying from point $A$ to point $B$.
In an ideal case, the trajectory of its flight
would be a straight line. However, for the reasons
which are mostly unknown to us, we observe the deviations
which we (due to lack of knowledge, or due to lack of
a full and precise model of the bumblebee's motion)
consider as random.

\begin{figure*}
\centering
\includegraphics[width=0.6\textwidth,angle=0]{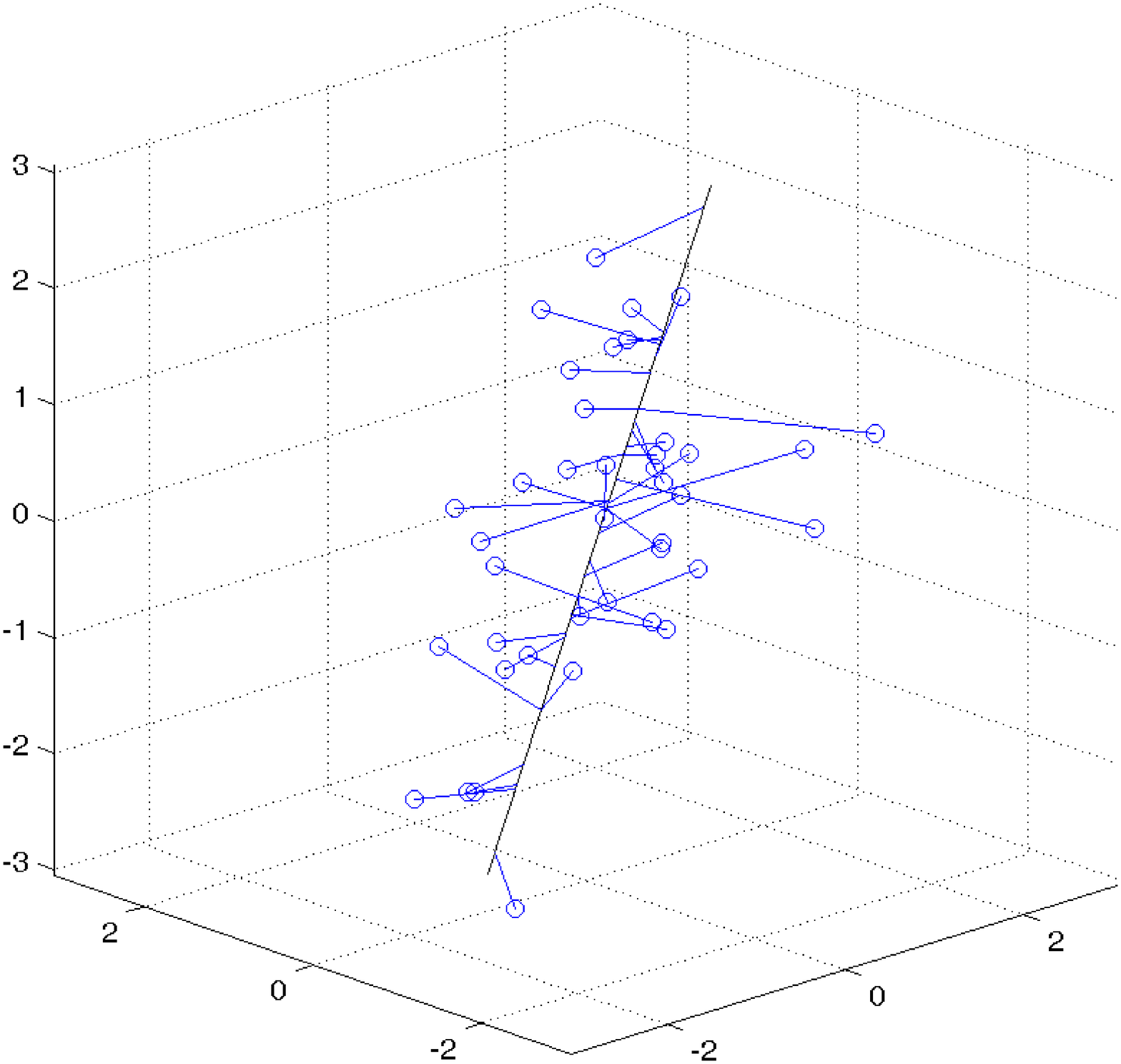}
\caption{The Flight of the Bumblebee.
(Inspired by a figure from \cite{princomp-demo}).
}\label{fig:bee}
\end{figure*}

Having a set of positions of the bumblebee (depicted
as discrete points in Fig.~\ref{fig:bee}), we can
try to make a conclusion about the main direction
of its flight. The simplest model is a straight line.
The only way to determine the parameters of the
equations of a fitting straight line in 3D space is
using the method of orthogonal regression (total
least squares). The result is shown in Fig.~\ref{fig:bee}
as a straight line.

Of course, the general trajectory of a bumblebee
can be more complex than a straight line.
Although in the case of orthogonal regression
we do not need to postulate the roles of variables,
we still need to postulate the structure of the model --
we have to decide which line is more ``promising''
for modeling the bumblebee's flight -- it could be
a straight line, or a conic section, or some spiral,
or something else.

The ``bumblebee flight'' problem can be considered
also as a problem of determining the direction
of the flight of an unknown aircraft, or determining
the state of economy at some time instant in the
future based on the observations in the past.

\section{The case of the Central European countries of the V4 group}

\begin{figure}
\centering
\noindent
  \subfigure[GDP of the V4 countries.]{\label{fig:GDPv4}
		\includegraphics[width=0.3\textwidth,angle=0]{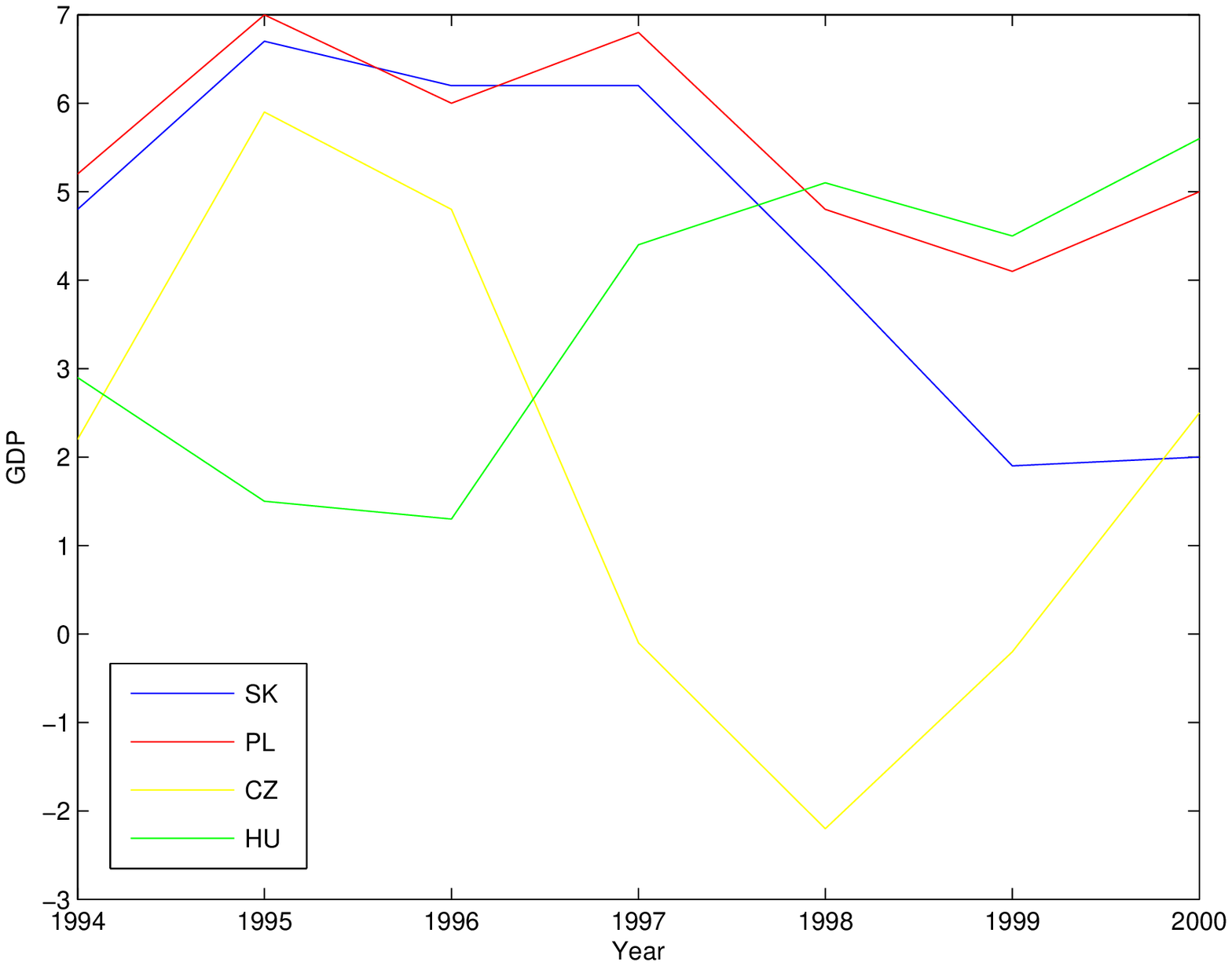}}
  \subfigure[Inflation of the V4 countries.]{\label{fig:InflationV4}
                \includegraphics[width=0.3\textwidth,angle=0]{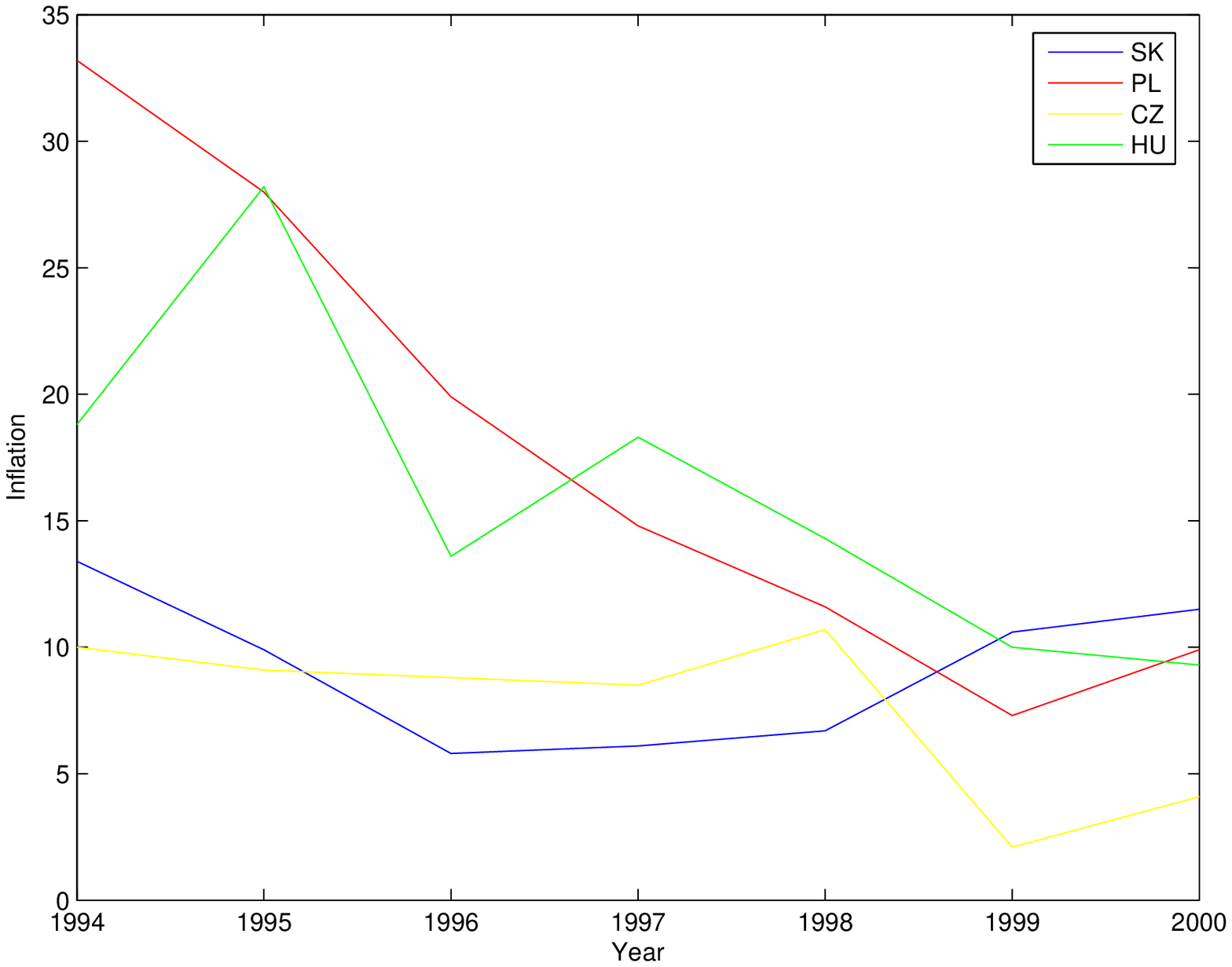}}
  \subfigure[Unemployment in the V4 countries.]{\label{fig:UnemploymentV4}
                \includegraphics[width=0.3\textwidth,angle=0]{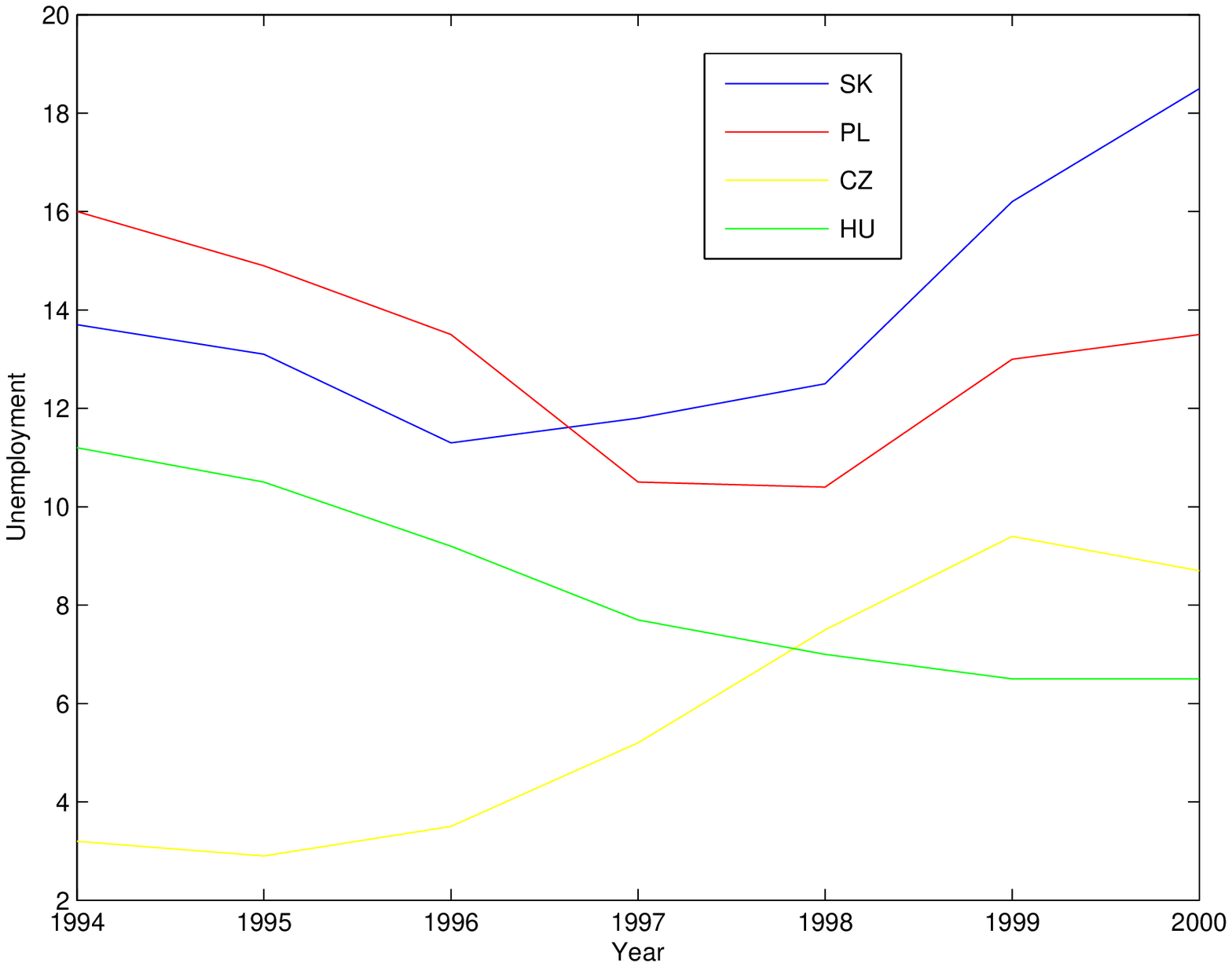}}
  \caption{Classical 2D visualization of economic indicators of the V4 countries.}
\end{figure}

The traditional approach to visualization of the change
of GDP, inflation, and unemployment uses 2D plots --
see Fig.~\ref{fig:GDPv4}, Fig.~\ref{fig:InflationV4}, and Fig.~\ref{fig:UnemploymentV4}.
Although such visualization provides partial information
about each particular aspect of a particular economy,
and even allows comparisons of a given particular
indicator for several economies, it does not provide
any information about the economy as a whole.

This is why we plot the data for each country of the V4 group --
Slovak Republic (SK), Czech Republic (CZ), Hungary (HU), and Poland (PL) --
in the state space with three coordinates:
GDP, inflation, and unemployment.
We use the data for the period from 1994 to 2000
from the MESA~10 report \cite{Jakoby-MESA10}.

\begin{figure}
\centering
\par\vspace*{-1cm}\par
\noindent
  \subfigure[Czech Republic]{\label{fig:CZ-phase-trajectory}
		\includegraphics[width=0.35\textwidth,angle=0]{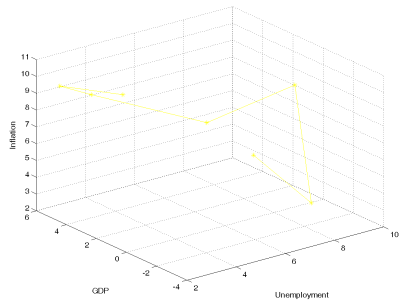}}
  \subfigure[Hungary]{\label{fig:HU-phase-trajectory}
		\includegraphics[width=0.35\textwidth,angle=0]{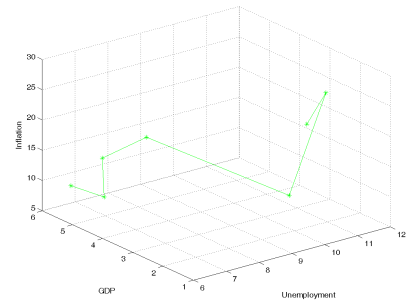}}
  \subfigure[Poland]{\label{fig:PL-phase-trajectory}
		\includegraphics[width=0.35\textwidth,angle=0]{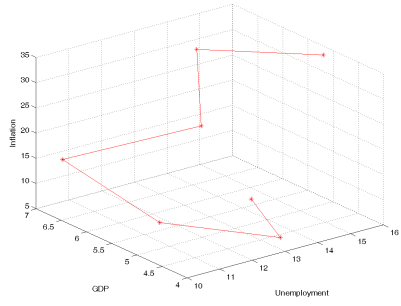}}
  \subfigure[Slovak Republic]{\label{fig:SK-phase-trajectory}
		\includegraphics[width=0.35\textwidth,angle=0]{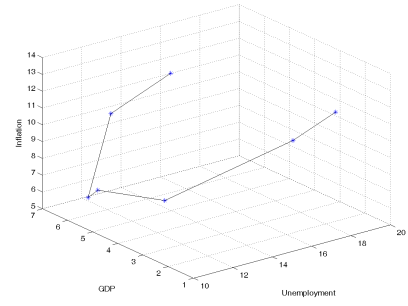}}
  \caption{Phase trajectories of the national economies of the V4 countries in state space.}
\end{figure}

One can see that, in conrast with the traditional 2D plots,
the phase trajectories shown in
Figs.~\ref{fig:CZ-phase-trajectory}--\ref{fig:SK-phase-trajectory}
nicely show how the national economies developed in time.

The most interesting observation regarding these
phase trajectories is that for each particular
country
its phase trajectory lies approximately in one plane.
This observation indicates that we can associated
such planes in state space with particular economies,
and that the global properties of each particular national
economy as a whole are described by the associated plane
(or, in other word, by its normal vector and
the data centroid, which is a point belonging to the plane).

The values of normal vectors and centroids of the planes
describing national economies of the V4 countries
are listed in Table~\ref{tab:NCE}. As it is well known
from the analytic geometry, each plane is uniquely described
by these data. In the last column the values of the
total error of approximation are given.

It is worth mentioning that the fitting of the state space
data points to planes was done using total least squares (TLS)
method, implemented in the form of the principal component
analysis \cite{Schuermans-equiv-TLSM_PCA-2005},
although in many cases the SVD decomposition approach
is used \cite{Golub-van-Loan-1980,Nievergelt-review-1994}.
 The description of our function for MATLAB
is given in the Appendix.

\begin{table}
\caption{Economic indicators of the V4 countries, 1994--2000}
\centering
{\footnotesize
\subtable[Unemployment in the V4 countries -- in percents]{
\begin{tabular}{llllllll}
\hline
~	& 1994	& 1995	& 1996	& 1997	& 1998	& 1999	& 2000 \\
\hline
CZ	& 3.2	& 2.9	& 3.5	& 5.2	& 7.5	& 9.4	& 8.7 \\
HU	& 11.2	& 10.5	& 9.2	& 7.7	& 7	& 6.5	& 6.5 \\
PL	& 16.0	& 14.9	& 13.5	& 10.5	& 10.4	& 13.0	& 13.5 \\
SK	& 13.7	& 13.1	& 11.3	& 11.8	& 12.5	& 16.2	& 18.5 \\
\hline
\end{tabular}
}
\subtable[GDP of the V4 countries -- change in percents]{
\begin{tabular}{llllllll}
\hline
	& 1994	& 1995	& 1996	& 1997	& 1998	& 1999	& 2000\\
\hline
CZ	& 2.2	& 5.9	& 4.8	& -0.1	& -2.2	& -0.2	& 2.5 \\
HU	& 2.9	& 1.5	& 1.3	& 4.4	& 5.1	& 4.5	& 5.6 \\
PL	& 5.2	& 7	& 6.0	& 6.8	& 4.8	& 4.1	& 5.0 \\
SK	& 4.8	& 6.7	& 6.2	& 6.2	& 4.1	& 1.9	& 2.0 \\
\hline
\end{tabular}
}
\subtable[Inflation in the V4 countries -- in percents]{
\begin{tabular}{llllllll}
\hline
~	& 1994	& 1995	& 1996	& 1997	& 1998	& 1999	& 2000 \\
\hline
CZ	& 10	& 9.1	& 8.8	& 8.5	& 10.7	& 2.1	& 4.1 \\
HU	& 18.8	& 28.2	& 23.6	& 18.3	& 14.3	& 10	& 9.3 \\
PL	& 33.2	& 28	& 19.9	& 14.8	& 11.6	& 7.3	& 9.9 \\
SK	& 13.4	& 9.9	& 5.8	& 6.1	& 6.7	& 10.6	& 11.5 \\
\hline
\end{tabular}
}
}
\end{table}

%
%
%
%


%
%


%
%
%


\begin{figure}
\centering
\par\vspace*{-1cm}\par
\noindent
  \subfigure[Czech Republic]{\label{fig:CZ3D}
		\includegraphics[width=0.35\textwidth,angle=0]{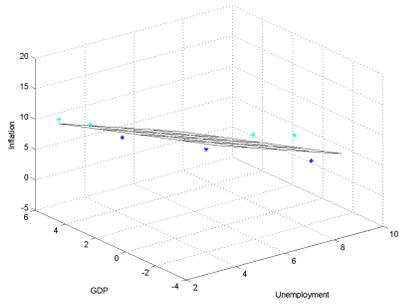}}
  \subfigure[Hungary]{\label{fig:HU3D}
		\includegraphics[width=0.35\textwidth,angle=0]{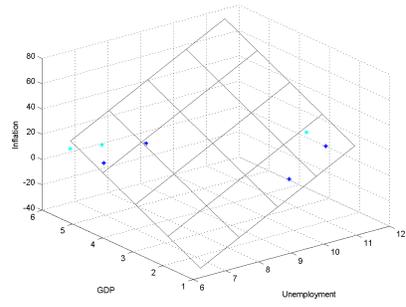}}
  \subfigure[Poland]{\label{fig:PL3D}
		\includegraphics[width=0.35\textwidth,angle=0]{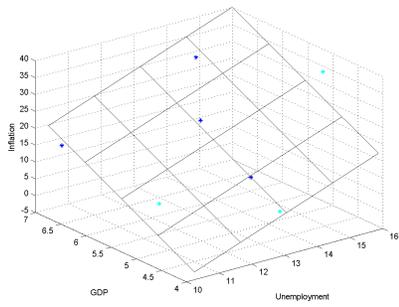}}
  \subfigure[Slovak Republic]{\label{fig:SK3D}
		\includegraphics[width=0.35\textwidth,angle=0]{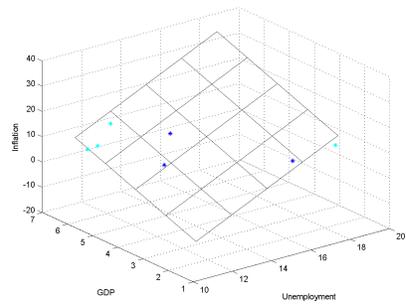}}
  \caption{Planes of the national economies of the V4 countries in state space.}
\end{figure}

%
%
%

\begin{table}[h!]
\caption{Normal vectors, centroids, and errors for the planes
corresponding to the countries of the V4 group.}\label{tab:NCE}
\par\vspace*{1.5ex}\par
\centering
\begin{tabular}{|c|c|c|c|}
\hline
\rule{0mm}{3ex}\textbf{Country}	& \textbf{Normal vector} & \textbf{Centroid} &	\textbf{Error} \\[1.5ex]
\hline
\rule{0mm}{3ex}
SK	& (0.6704,~0.7195,~-0.1811)	& (13.8714,~4.5571,~9.1429)	& 4.2633 \\[1.5ex]
PL		& (-0.4083,~-0.9059,~0.1123)	& (13.1143,~5.5571,~17.8143)	& 4.3106 \\[1.5ex]
CZ	& (0.7632,~0.4525,~0.4612)	& (5.7714,~1.8429,~7.6143)	& 4.6111 \\[1.5ex]
HU		& (0.7362,~0.6745,~-0.0545)	& (8.3714,~3.6143,~16.0714)	& 3.7431 \\[1.5ex]
\hline
\end{tabular}
\end{table}

\section{Conclusions}

In this article we presented a new approach to description of
national economies. This approach is based on the state space viewpoint,
which is used mostly in the theory of dynamical systems
and in the control theory. Gross domestic product,
inflation, and unemployment rates were taken as state
variables. We demonstrated that for the considered
period of time the phase trajectory of each of the V4 countries
(Slovak Republic, Czech Republic, Hungary, and Poland)
lies approximately in one plane, so that the economic
development of each country can be assocated with
a corresponding plane in the state space.

The suggested approach opens a way to a new set
of economic indicators. Among possible indicators
of economies we can mention, for example, normal vectors
of national economies, various plane slopes,
2D angles between the planes corresponding to
different economies, etc.

The tool used for computations is orthogonal regression
(alias orthogonal distance regression,
alias total least squares method),
and we also gave general arguments for using orthogonal regression
instead of the classical regression
based on the least squares method.

\section*{Acknowledgements}

This work was partially supported by grants
VEGA 1/2165/05, VEGA 1/3132/06, VEGA 1/2179/05,
VEGA 1/2160/05, and grant number 49s3
under Aktion Osterreich-Slowakei.
The authors are grateful to Prof. Paul O'Leary,
Montanuniversit\"at Leoben, Austria,
for his advices and provided materials.

\appendix

\section{Brief description of the MATLAB routine}

We created a MATLAB routine called \verb+fit_3D_data+
for performing computations used in our article.
Since our routine can be used in many other
investigations utilizing the orthogonal regression
(total least squares) method, we published it
at MATLAB File Exchange. The URI for download is:\\
\verb+     http://www.mathworks.com/matlabcentral/fileexchange/+ \\
\verb+              loadFile.do?objectId=12395&objectType=file+ \\
The description of \verb+fit_3D_data+
and its parameters is given below.

\begin{verbatim}

function [Err,N,P] = fit_3D_data(XData, YData, ZData,
                                 geometry, visualization, sod)
%
% [Err, N, P] = fit_3D_data(XData, YData, ZData,
%                           geometry, visualization, sod)
%
% Orthogonal Linear Regression in 3D-space
% by using Principal Components Analysis
%
% This is a wrapper function to some pieces of the code from
% the Statistics Toolbox demo titled "Fitting an Orthogonal
% Regression Using Principal Components Analysis"
% (http:/ /www.mathworks.com/products/statistics/
%  demos.html?file=/products/demos/shipping/stats/orthoregdemo.html),
% which is Copyright by the MathWorks, Inc.
%

% Input parameters:
%  - XData: input data block -- x: axis
%  - YData: input data block -- y: axis
%  - ZData: input data block -- z: axis
%  - geometry: type of approximation ('line','plane')
%  - visualization: figure ('on','off') -- default is 'on'
%  - sod: show orthogonal distances ('on','off') -- default is 'on'
%
% Return parameters:
%  - Err: error of approximation - sum of orthogonal distances
%  - N: normal vector for plane, direction vector for line
%  - P: point on plane or line in 3D space
%
% Example:
%
% >> XD = [4.8 6.7 6.2 6.2 4.1 1.9 2.0]';
% >> YD = [13.4 9.9 5.8 6.1 6.7 10.6 11.5]';
% >> ZD = [13.7 13.1 11.3 11.8 12.5 16.2 18.5]';
% >> fit_3D_data(XD,YD,ZD,'line','on','on');
% >> fit_3D_data(XD,YD,ZD,'plane','on','off');
%
% Note: Written for Matlab 7.0 (R14) with Statistics Toolbox
%
% We sincerely thank Peter Perkins, the author of the demo,
% and John D'Errico for their comments.
%
% Ivo Petras, Igor Podlubny, May 2006
% (ivo.petras@tuke.sk, igor.podlubny@tuke.sk)
\end{verbatim}

\end{document}